\documentclass[leqno,11pt]{article}
\usepackage{amsmath, amssymb, amsthm}
\usepackage[all]{xy}

\usepackage{slashed}

\newtheorem{theorem}{Theorem}[section]

\newtheorem{proposition}[theorem]{Proposition}
\newtheorem{corollary}[theorem]{Corollary}
\newtheorem{definition}[theorem]{Definition\rm}

\newcommand{\ggot}{\ensuremath{\mathfrak{g}}}

\newcommand{\sogot}{\ensuremath{\mathfrak{so}}}

\newcommand{\Acal}{\ensuremath{\mathcal{A}}}

\newcommand{\Ccal}{\ensuremath{\mathcal{C}}}

\newcommand{\Ecal}{\ensuremath{\mathcal{E}}}

\newcommand{\Kcal}{\ensuremath{\mathcal{K}}}

\newcommand{\Qcal}{\ensuremath{\mathcal{Q}}}

\newcommand{\Vcal}{\ensuremath{\mathcal{V}}}

\newcommand{\C}{\ensuremath{\mathbb{C}}}

\newcommand{\R}{\ensuremath{\mathbb{R}}}

\newcommand{\f}{\ensuremath{\mathcal{C}^{\infty}}}
\newcommand{\fgene}{\ensuremath{\mathcal{C}^{-\infty}}}

\newcommand{\indice}{\ensuremath{\hbox{\rm Indice}}}
\newcommand{\ch}{\ensuremath{\hbox{\rm Ch}}}

\newcommand{\Cl}{\ensuremath{\hbox{\rm Cl}}}

\newcommand{\T}{\ensuremath{\hbox{\bf T}}}
\newcommand{\tr}{\ensuremath{\hbox{\bf Tr}}}

\newcommand{\K}{\ensuremath{\hbox{\bf K}}}

\newcommand{\End}{\ensuremath{\hbox{\rm End}}}

\def \so {{\rm SO}}
\def \spin {{\rm Spin}}
\def \su {{\rm SU}}
\def \pu {{\rm PU}}

\def \wG {{\widetilde{G}}}

\def \wgamma {{\widehat{\Gamma}}}

\def \wg {{\widetilde{g}}}

\def \wggot {{\widetilde{\mathfrak{g}}}}

\begin{document}

\title{Index of projective elliptic operators}

\author{Paul-Emile PARADAN\footnote{Institut Montpelli\'erain Alexander Grothendieck, CNRS UMR 5149,
Universit\'e de Montpellier, \texttt{paul-emile.paradan@umontpellier.fr}}}

\maketitle

\begin{abstract}
Mathai, Melrose, and Singer introduced the notion of projective elliptic operators on manifolds 
equipped with an Azumaya bundle. In this note we compute the equivariant index of transversally 
elliptic operators that are the pullback of projective elliptic operators on the trivialization of the 
Azumaya bundle. It encompasses the fractional index formula of projective elliptic operator by 
Mathai-Melrose-Singer.
\end{abstract}

\section{Introduction}

In \cite{MMS05,MMS06,MMS08,MMS09}, R. Melrose, V. Mathai and I. M. Singer studied the notion of projective elliptic operators and the corresponding index problem. There were interested to the situation were a compact manifold $M$ carries an Azumaya bundle $\Acal\to M$ of rank $N\times N$. In this setting we have a $\pu_N$-principal bundle $P_\Acal\to M$ which corresponds to the trivialization of $\Acal$, 
$$
\Acal\simeq P_\Acal\times_{\pu_N}\End(\C^N),
$$
and one works with the central extension $1\to \Gamma_N \to \su_N\to \pu_N\to 1$.

An $\Acal$-projective elliptic operator on $M$ can be understood as an elliptic operator on the orbifold $P_\Acal/\su_N\simeq M$. Thanks to the work of Atiyah-Singer \cite{Ati74} and Kawasaki \cite{K81},  a good way to think to an elliptic operator $D$ on the orbifold $P_\Acal/\su_N$ is to view it, 
after pulling back by the projection $P_\Acal\to M$, as a $\su_N$-transversally elliptic operators $\tilde D$ on $P_\Acal$. 

In \cite{MMS06}, the authors define the analytical index, $\indice_a(D)$, of an $\Acal$-projective elliptic operator $D$ on $M$, and they prove that one has an Atiyah-Singer type formula for this invariant which shows in particular that $\indice_a(D)$ is a rational number.
%

In \cite{MMS08},  the authors explained the link between the analytical index of $D$ and the equivariant index  of the transversally elliptic operator 
$\tilde{D}$, denoted $\indice_{\su_N}(\tilde{D})$. Recall that $\indice_{\su_N}(\tilde{D})$ can be viewed as an invariant distribution on $\su_N$ and Atiyah-Singer theory tell us that $\indice_{\su_N}(\tilde{D})$ is supported on the subgroup $\Gamma_N\subset \su_N$. The main result of \cite{MMS08} is that 
\begin{equation}\label{eq:intro-2}
\indice_a(D)=\left\langle \indice_{\su_N}(\tilde{D}),\varphi\right\rangle_{\su_N}
\end{equation}
when $\varphi$ is a smooth function on $\su_N$ which is constantly equal to $1$ around the identity element and has small enough support.

The above formula gives ``the coefficient of the Dirac function'' in the distribution $\indice_{\su_N}(\tilde{D})$. 
The main purpose of the present note is to compute the ``higher order derivatives of the Dirac function'' in the distribution $\indice_{\su_N}(\tilde{D})$ (see  Theorem \ref{theo-intro}). Our computation  generalizes the previous result obtained by Yamashita \cite{Y13} in the case of the projective Dirac operator. 

\section{Index of transversally elliptic operators}

We consider the following setting :

\begin{enumerate}
\item a $G$-principal bundle $\pi: P\to M$, where $G$ is a compact Lie group,
\item a central extension $1\to \Gamma \to \wG\to G\to 1$ where $\Gamma$ is finite.
\end{enumerate}

We are interested to the $\wG$-equivariant pseudo-differential operators on $P$ that are 
transversally elliptic \cite{Ati74}. We are here in a very particular situation: the sugbroup 
$\Gamma$ acts trivially on $P$ and  the group $G$ acts freely on $P$, so the 
set $\T^*_\wG P=\T^*_G P$ formed by the covectors orthogonal to the orbits is 
a sub-bundle of $\T^*P$, and the quotient by $G$ defines a principal bundle $\T^*_G P\to \T^* M$.

The index of a $\wG$-transversally elliptic operator on $P$, says $\tilde{D}$, depends uniquely of the class defined by its principal symbol $\sigma(\tilde{D})$ in the group $\K^0_{\wG}(\T^*_G P)$ (see \cite{Ati74}). Hence the analytical index defines a morphism
$$
\indice_{\wG}: \K^0_{\wG}(\T^*_{G}P)\longrightarrow \fgene(\wG)^{\rm Ad},
$$
and our aim is to compute it by cohomological formulas. Our main tool will be delocalized 
formulas obtained by Berline-Paradan-Vergne in \cite{BV1,BV2,PV}.

A general property of the equivariant index states that for any class $\sigma \in  \K^0_{\wG}(\T^*_{G}P)$ 
the distribution $\indice_{\wG}(\sigma)$ is supported on elements $\wg\in \wG$ such that 
$P^{\tilde g}\neq \emptyset$. In our case, such elements belong to the finite subgroup $\Gamma$, 
so we have a decomposition
$$
\indice_{\wG}(\sigma)=\sum_{\gamma\in\Gamma} \Qcal_\gamma(\sigma)
$$
where $\Qcal_\gamma(\sigma)$ is an Ad-invariant distribution on $\wG$ supported at the central element $\gamma$. Recall that the envelopping algebra $U(\ggot)=U(\wggot)$ is canonically identified with the algebra of distributions on $\wG$ supported at the identity. So the center $Z(\ggot)=U(\ggot)^\ggot$ corresponds to the Ad-invariant distributions on $\wG$ supported at the identity.

Let $\delta_\gamma\in \fgene(\wG)^{\rm Ad}$ be the Dirac distribution supported at the central element $\gamma\in \Gamma$. The convolution map 
$$
T\in Z(\ggot)\longmapsto T\star \delta_\gamma
$$
defines a bijection between $Z(\ggot)$ and the {\rm Ad}-invariant distributions on $\wG$ supported at $\gamma$. As the invariant distribution $\Qcal_\gamma(\sigma)$ is supported at $\gamma$, there exists an element $T_{\gamma}(\sigma)\in Z(\ggot)$ such that $\Qcal_\gamma(\sigma)=T_{\gamma}(\sigma)\star  \delta_\gamma$.

The exponential  map $\exp:\ggot\to G$ 
defines a linear isomorphism\footnote{It it not a morphism of algebras !} 
$$
\exp_*: S(\ggot)^G\to Z(\ggot)
$$
 where $S(\ggot)^G$ is viewed as the algebra of Ad-invariant distributions on $\ggot$ supported at $0$.  
 Our main purpose is to give a cohomological formula for the element
$\exp_*^{-1}\left(T_{\gamma}(\sigma)\right)\in S(\ggot)^G$.

\section{Chern-Weil morphism and twisted Chern characters}

Let $E_1,\cdots,E_r$ be a basis of $\ggot$,  and let $\theta=\sum_{i=1}^r \theta_i\otimes E_i\in (\Acal^1(P)\otimes \ggot)^G$ be a connection on the principal bundle $P\to M$. 
Its curvature $\Theta=\sum_{i=1}^r \Theta_i\otimes E_i$ belongs to $(\Acal_{hor}^{+}(P)\otimes \ggot)^G$, where $\Acal_{hor}^{+}(P)$ is the algebra of horizontal 
forms of even degree on $P$.  A polynomial function $P(X)=P(X_1,\ldots, X_r)$ on $\ggot$ can be 
evaluated at $\Theta$: the corresponding element $P(\Theta)=P(\Theta_1,\ldots,\Theta_r)$ is a 
closed differential form of $\Acal_{bas}^{+}(P)\simeq \Acal^{+}(M)$ when $P$ is $G$-invariant. We can define
$\varphi(\Theta)$ for any smooth function $\varphi:\ggot\to \C$, by using its Taylor series at $0$. If we denote $H_{dR}(M)$ the de Rham cohomology of $M$, the map 
$\varphi\in \f(\ggot)^{\rm Ad}\mapsto \varphi(\Theta)\in H_{dR}(M)$ is the Chern-Weil morphism (it is independent of the choice of the connection).

Another way to look at the Chern-Weil morphism is to consider the exponential $e^\Theta\in\left(\Acal_{hor}^{+}(P)\otimes S(\ggot)\right)^G$. If we view the polynomial algebra $S(\ggot)$ as the algebra of distributions on $\ggot$ supported at $0$, we see that 
$$\varphi(\Theta):=\left\langle e^\Theta,\varphi\right\rangle_\ggot.$$
Note that $\varphi(\Theta)=1$ if $\varphi$ is equal to $1$ in a neighborhood of $0$.

\begin{definition}
For any closed form $\alpha$ on $\T^*M$ with compact support, the expression
$\int_{{\rm T}^*M}\alpha\wedge e^\Theta$ defines an element of $S(\ggot)^G$ through the relation
$$
\left\langle\int_{{\rm T}^*M}\alpha\wedge e^\Theta,\varphi\right\rangle_\ggot:=\int_{{\rm T}^*M}\alpha\wedge \overline{\varphi}(\Theta),
$$
that holds for any $\varphi\in\f(\ggot)$. Here $\overline{\varphi}=\int_G g\cdot\varphi \,dg$ is the average of $\varphi$ relatively to the Haar mesure on $G$ of volume $1$, and the 
 forms $ \overline{\varphi}(\Theta)\in \Acal^+(M)$ are also considered as forms on $\T^* M$. Note that $\int_{{\rm T}^*M}\alpha\wedge e^\Theta\in S(\ggot)^G$ is independent of the choice of the connection.
\end{definition}

\medskip

Now we explain how are constructed the twisted Chern characters $\ch_{\gamma}(\sigma)$ 
of a $\wG$-transversally elliptic symbol on $P$.

Let   $\sigma\in \Ccal^{\infty}(\T^*P,\hom(p^*\Ecal^+,p^*\Ecal^-))$ be a $\wG$-transversally elliptic symbol: 
here $\Ecal^\pm$ are $\wG$-complex vector bundles on $P$ and $p:\T^*P\to P$ is the projection.  By definition the set 
$$
\Kcal_\sigma:= \left\{(x,\xi)\in \T^*_G P\ \vert\  \sigma(x,\xi):\Ecal^{-}_{x}\to\Ecal^{+}_{x}\ \mathrm{is\ not\ invertible}\right\}
$$
is compact. 

We start with a $\wG$-equivariant connection $\nabla^{\Ecal^+}$ on the vector bundle $\Ecal^+\to P$. 
The pull-back $\nabla^{p^*\Ecal^+}:=p^*\nabla^{\Ecal^+}$ is then a connection on $p^*\Ecal^+$ viewed as a vector bundle on the manifold 
$\T^*_G P$. Since $\Kcal_\sigma$ is compact, we can define on the vector bundle $p^*\Ecal^-\to \T^*_G P$ a connection $\nabla^{p^*\Ecal^-}$ such that the following relation 
\begin{equation}\label{eq:nabla}
\nabla^{p^*\Ecal^-}=\sigma\circ \nabla^{p^*\Ecal^+}\circ \sigma^{-1}
\end{equation}
holds outside a compact subset of $\T^*_G P$.

Let $R^+(X), R^-(X), X\in\ggot$ be the equivariant curvatures of the  connections $\nabla^{p^*\Ecal^+}$ and $\nabla^{p^*\Ecal^-}$ (see \cite{BGV}, section 7). We consider the equivariant Chern character, twisted by the central element $\gamma\in \Gamma$:
$$
\ch^{\wG}_\gamma(\sigma)(X):= \tr\left(\gamma^{\Ecal^+}e^{R^+(X)}\right)- \tr\left(\gamma^{\Ecal^-}e^{R^-(X)}\right),\quad X\in\ggot.
$$
In this formula, $\gamma^{\Ecal^\pm}$ denotes the linear action of $\gamma$ on the fibers of the bundles $\Ecal^\pm$.
Thanks to (\ref{eq:nabla}), the closed equivariant form $\ch^{\wG}_\gamma(\sigma)$ has a compact support on $\T^*_G P$. It defines a class, still denoted 
$\ch^{\wG}_\gamma(\sigma)$, in the equivariant cohomology group with compact support $H^{\infty}_{\wG, c}(\T^*_G P)$. Since the finite subgroup $\Gamma$ acts trivially on $P$, we have a canonical isomorphism between $H^{\infty}_{\wG, c}(\T^*_G P)$ and $H^{\infty}_{G, c}(\T^*_G P)$.

\begin{definition}\label{def:chern}
Let $H_{dR,c}(\T^*M)$ be the de Rham cohomology of $\T^*M$ with compact support. The twisted Chern characters 
$\ch_{\gamma}(\sigma)\in H_{dR,c}(\T^*M)$ is defined as the image of $\ch^{\wG}_\gamma(\sigma)$ under the Chern-Weil isomorphism 
$$
 H^{\infty}_{G, c}(\T^*_G P)\to H_{dR,c}(\T^*M)
$$ 
that is associated to the principal $G$-bundle $\T^*_G P\to \T^*M$.
\end{definition}

\section{Main result}

Let $R_M$ is the curvature of the tangent bundle $\T M\to M$. The $\widehat{\mathrm{A}}$-class of $M$ is normalized as follows
$$
\widehat{\mathrm{A}}(M)=\det{}^{1/2}\left(\frac{R_M}{e^\frac{R_M}{2}-e^\frac{-R_M}{2}}\right) \  \in\  H_{dR}(M).
$$

The following Theorem is the main result of this note.

\begin{theorem}\label{theo-intro} 
Let $\sigma\in\K^0_{\wG}(\T^*_{G}P)$. We have 
$\indice_{\wG}(\sigma)=\sum_{\gamma\in\Gamma} T_{\gamma}(\sigma)\star \delta_\gamma$ where
$$
T_{\gamma}(\sigma)=(2i\pi)^{-{\rm dim} M}\exp_*\left(\int_{{\rm T}^*M}\widehat{\mathrm{A}}(M)^2
\wedge\ch_{\gamma}(\sigma) \wedge e^\Theta\right).
$$
Here  $\ch_{\gamma}(\sigma)$ is the twisted Chern character (see Definition \ref{def:chern}).
\end{theorem}

\begin{corollary} We have
$$
\left\langle\indice_{\wG}(\sigma),\varphi\right\rangle_{\wG}=(2i\pi)^{-{\rm dim} M} \int_{{\rm T}^*M}\widehat{\mathrm{A}}(M)^2
\wedge\ch_{\gamma}(\sigma)
$$
if $\varphi\in\f(\wG)$ is a fonction constantly equal to $1$ around $\gamma$ and has small enough support.

\end{corollary}

The proof of Theorem \ref{theo-intro} follows from a direct application of the delocalized formulas obtained by 
Berline-Paradan-Vergne in \cite{BV1,BV2,PV}. Here we just explain how one can obtain easily this result 
when $\wG=G$. In this setting the map 
$\pi:P\to M$ defines an isomorphism $\sigma\mapsto \pi^*\sigma, \K^0(\T^*M)\to \K^0_{G}(\T^*_{G}P)$ at the level of 
$\K$-groups, and the distributions $\indice_{G}(\pi^*\sigma)$ are supported only at the identity element.

 For any irreducible representation $V_\lambda$ (parametrized by $\lambda\in\hat{G}$) we denote 
$\chi_\lambda$ its character and $\Vcal_\lambda$ the complex vector bundles 
$P\times_G V_\lambda\to M$. Theorem 3.1 of Atiyah in \cite{Ati74}  tells us that
$$
\langle \indice_G(\pi^*\sigma), \chi_\lambda\rangle_G=
\indice(\sigma\otimes \Vcal_\lambda)\nonumber\\ 
=(2i\pi)^{-{\rm dim} M}\int_{{\rm T}^* M} \hat{A}(M)^2\wedge \ch(\sigma)\wedge\ch(\Vcal_\lambda).
$$
The Chern character of the vector bundle $\Vcal_\lambda$ is defined by the Chern-Weil morphism:
$\ch(\Vcal_\lambda)= \langle e^{\Theta}, \exp^*(\chi_\lambda)\rangle_\ggot$. Then
$$
\left\langle \indice_G(\pi^*\sigma), \chi_\lambda\right\rangle_G=(2i\pi)^{-{\rm dim} M}
\left\langle\int_{{\rm T}^* M}\hat{A}(M)^2\wedge \ch(\sigma)\wedge e^{\Theta},\exp^*(\chi_\lambda) \right\rangle_\ggot
$$
for any $\lambda\in\hat{G}$. We can conclude that
$$
\indice_G(\pi^*\sigma)= (2i\pi)^{-{\rm dim} M}\exp_*\left(\int_{{\rm T}^*M}\widehat{\mathrm{A}}(M)^2
\wedge\ch(\sigma) \wedge e^\Theta\right)\star \delta_e.
$$

\medskip

If we consider the group $\wgamma$ of characters of the finite abelian group $\Gamma$, we can decompose a 
$\wG$-transversally elliptic symbol $\sigma\in\Ccal^{\infty}(\T^*P,\hom(p^*\Ecal^+,p^*\Ecal^-))$ as 
$\sigma=\oplus_{\chi\in\wgamma}\ \sigma_\chi$,  
where $\sigma_\chi\in\Ccal^{\infty}(\T^*P,\hom(p^*\Ecal^+_\chi,p^*\Ecal^-_\chi))$ is a $\wG$-transversally elliptic symbol on $P$. Here we take $\Ecal^\pm_\chi$ as the subbundle of $\Ecal^\pm$ where $\Gamma$ acts trough the character $\chi$.  From Definition \ref{def:chern}, it is 
obvious that the twisted Chern character $\ch_\gamma(\sigma)$ admits the decomposition
$$
\ch_\gamma(\sigma)=\sum_{\chi\in\wgamma}\, \langle\chi,\gamma\rangle\, \ch_e(\sigma_\chi),
$$ 
where $\langle-,-\rangle: \wgamma\times \Gamma\to \C$ is the duality bracket.

We have a reformulation of Theorem \ref{theo-intro}.

\begin{theorem}\label{theo-intro-bis} 
Let $\sigma\in\K^0_{\wG}(\T^*_{G}P)$ with decomposition $\sigma=\oplus_{\chi\in\wgamma}\ \sigma_\chi$. We have 
$$
\indice_{\wG}(\sigma)=\sum_{(\chi,\gamma)\in\wgamma\times\Gamma} \, \langle\chi,\gamma\rangle\,T_{e}(\sigma_\chi)\star \delta_\gamma
$$
where
$T_{e}(\sigma_\chi)=(2i\pi)^{-{\rm dim} M}\exp_*\left(\int_{{\rm T}^*M}\widehat{\mathrm{A}}(M)^2
\wedge\ch_{e}(\sigma_\chi) \wedge e^\Theta\right)$.

\end{theorem}

\section{Projective elliptic operators}

Let us come back to the setting of an  $\Acal$-projective elliptic operator $D$ on $M$. Let $\tilde D$ be its pullback on $P_\Acal$: it is a $\su_N$-transversally elliptic operator, and we denote $\sigma(\tilde{D})\in \K^0_{\su_N}(\T^*_{\pu_N}P_\Acal)$ its principal symbol. Here the action of the center $\Gamma_N\subset \su_N$ on 
$\sigma(\tilde{D})$ is prescribed by the canonical character $\chi_o : \Gamma_N\to \C$. In other words,  
$\sigma(\tilde{D})=\sigma(\tilde{D})_{\chi_o}$.

In this case Theorem \ref{theo-intro-bis} gives

\begin{proposition}\label{prop-intro} Let $D$ be an $\Acal$-projective elliptic operator on $M$. Let $\tilde D$ be its pullback on $P_\Acal$. We have the following relation in 
 $\fgene(\su_N)^{\rm Ad}$,
$$
\indice_{\su_N}(\tilde{D})=T_D\star \left(\sum_{z\in \Gamma_N} \, \langle\chi_o,z\rangle \, \delta_z\right),
$$
where $T_D\in Z(\mathfrak{su}_N)$ is defined by the relation 
$$
T_D=(2i\pi)^{-{\rm dim} M}\exp_*\left(\int_{{\rm T}^*M}\widehat{\mathrm{A}}(M)^2\wedge\ch_e(\sigma(\tilde{D}))\wedge e^\Theta\right).
$$
Here $\Theta$ is the curvature of the $\pu_{N}$-principal bundle $P_\Acal\to M$.
\end{proposition}

\medskip

We finish this note by considering the particular case of the projective Dirac operator. 
Let $M$ be a compact Riemannian manifold of dimension $2n$ that is supposed oriented: here 
we work with the exemple of Azumaya bundle defined by the Clifford bundle ${\Acal}:=\Cl(\T M)_\C$.  
In this context  R. Melrose, V. Mathai and I. M. Singer showed that we have a natural 
projective Dirac operator $\slashed{\partial}_M^{\rm pr}$  on $M$ \cite{MMS06}.

In \cite{Y13}, Yamashita considers the pullback, denoted $\slashed{\partial}^{\rm su}_M$, of the projective Dirac operator $\slashed{\partial}_M^{\rm pr}$ relatively to the $\pu_N$-principal bundle $P_{\Acal}\to M$ (here $N=2^n$). The operator 
$\slashed{\partial}^{\rm su}_M$ is $\su_N$-transversally elliptic, and Yamashita computes the distribution 
on $\su_N$ defined by its equivariant index. The expression he obtains is similar to the one of 
Proposition \ref{prop-intro}, modulo a small mistake : the terms $\langle\chi_o,z\rangle$ are missing 
in his formula (see Corollary 6 in \cite{Y13}).

Instead of working with the $\pu_N$ principal bundle $P_{\Acal}\to M$, we can work with the $\so_{2n}$-principal bundle $P_{\so}\to M$ of oriented orthonormal frames, since it defines also a trivialization of $\Cl(\T M)_\C$. The pull-back of $\slashed{\partial}^{\rm pr}_M$ relatively to the projection $\pi:P_{\so}\to M$ is a $\spin_{2n}$-transversally elliptic operator on $P_{\so}$, that we denote $\slashed{\partial}^{\rm Spin}_M$. 

We recall the definition of its principal symbol $\sigma(\slashed{\partial}^{\rm Spin}_M)$. The kernel of 
$\T\pi:\T P\to \T M$ is the trivial sub-bundle isomorphic to $\sogot_{2n}\times P$. 
Let $\T^*_{\so}P\subset \T^*P$ be the orthogonal of $\sogot_{2n}\times P$. The bundle 
$\T^*_{\so}P$ admits 
an $\so_{2n}$-equivariant trivialisation $\alpha: \R^{2n}\times P\to \T^*_{\so}P$ defined as follows. 
For $f\in P$, the map $f:\R^{2n}\to \T_{\pi(f)}M$ is orthogonal and $(\T\pi\vert_f)^*: \T_{\pi(f)}^*M\to 
\T^*_{\so}P\vert_f$ is an isomorphism.  Hence the map 
$$
\alpha_f:=\T\pi^*\vert_f \circ (f^*)^{-1}: \R^{2n}\simeq ( \R^{2n})^*\longrightarrow \T^*_{\so}P\vert_f
$$
is a linear isomorphism. We see then that $P\times \R^{2n}\to \T^*_{\so}P, (f,x)\mapsto \alpha_f(x)$ 
is an $\so_{2n}$-equivariant trivialization.

Let $S_{2n}=S_{2n}^+ \oplus S_{2n}^-$ be the spinor representation. We denote 
$\mathrm{cl} : \R^{2n}\to \End(S_{2n})$ the clifford action (which is $\spin_{2n}$-equivariant). The symbol 
$\sigma(\slashed{\partial}^{\rm Spin}_M):\T^*_{\so}P\to\hom(S_{2n}^+, S_{2n}^-)$ is defined by the relations
$$
\sigma(\slashed{\partial}^{\rm Spin}_M)\vert_{(f,v)}
= \mathrm{cl}(\tilde{v}): S_{2n}^+ \to S_{2n}^-, \quad v\in  \T^*_{\so}P\vert_f,
$$
where $\tilde{v}=\alpha_f^{-1}(v)$.

\begin{proposition}\label{prop2-intro} We have the following equality  in $\fgene(\spin_{2n})^{\rm Ad}$:
$$
\indice_{\spin_{2n}}(\slashed{\partial}^{\rm Spin}_M)= T_M\star \delta_1 -  T_M\star \delta_{-1}
$$
where $T_M\in Z(\mathfrak{so}_{2n})$ is defined by the relation 
\begin{equation} \label{eq:T-M}
T_M=(2i\pi)^{-n}\exp_*\left(\int_{M}\widehat{\mathrm{A}}(M)\wedge e^\Theta\right).
\end{equation}
Here $\Theta$ is the curvature of the $\so_{2n}$-principal bundle $P_\so\to M$.
\end{proposition}

Proposition \ref{prop2-intro} follows from Theorem \ref{theo-intro-bis} as the Chern class 
$\ch(\slashed{\partial}^{\rm Spin}_M)$ satisfies the classical relation
$\ch_e(\sigma(\slashed{\partial}^{\rm Spin}_M))=(2i\pi)^{n}\,\widehat{\mathrm{A}}(M)^{-1}\wedge {\rm Thom}(\T^* M)$.

\end{document}